\documentclass[10pt]{article}
\usepackage{amssymb,amsbsy,amsmath,amsfonts,amssymb,amscd, mathrsfs}

\usepackage[english]{babel}
\usepackage[T1]{fontenc}
\usepackage{indentfirst}

\makeatletter
\@addtoreset{equation}{section}
\makeatother

\newtheorem{Theo}{Theorem}[section]
\newtheorem{Lem}[Theo]{Lemma}

\newcommand\N{\mathbb N}
\newcommand\R{\mathbb R}
\newcommand\T{\mathbb T}

\newcommand\Z{\mathbb Z}
\newcommand\E{{\mathbb E}\,}
\renewcommand\P{{\mathbb P}\,}
\newcommand\e{{\rm e}}

\newcommand\ind{{\rm 1\kern-.30em I}}

\begin{document}

\title{On some random thin sets of integers} 

\author{ 
Daniel \textsc{Li} -- Herv\'e \textsc{Queff\'elec} --  Luis \textsc{Rodr{\' \i}guez-Piazza} 
}

\date{}

\maketitle

\begin{abstract} 
\textit{We show how different random thin sets of integers may have different
behaviour. First, using a recent deviation inequality of Boucheron, Lugosi
and Massart, we give a simpler proof of one of our results in 
{\sl Some new thin sets of integers in Harmonic Analysis, Journal
d'Analyse Math\'ematique 86 (2002), 105--138}, namely that there exist 
$\frac{4}{3}$-Rider sets which are sets of uniform convergence and 
$\Lambda (q)$-sets for all $q < \infty $, but which are not Rosenthal sets.
In a second part, we show, using an older result of Kashin and Tzafriri that,
for $p > \frac{4}{3}$, the $p$-Rider sets which we had constructed in that
paper are almost surely not of uniform convergence.}
\end{abstract}

\vspace{2mm}
\selectlanguage{english}
\noindent \textbf{2000 MSC :}  \textit{primary : 43 A 46 ; secondary : 42 A 55 ; 42 A 61 }

\vspace{1mm}


\section{Introduction}
It is well-known that the Fourier series 
$S_n(f,x)  = \sum_{-n}^n \hat f(k) \e^{ikx}$ of a $2\pi$-periodic continuous function $f$ may be
badly behaved:  for example, it may diverge on a prescribed set of values of $x$ with measure zero. Similarly, the 
Fourier series of an integrable function may diverge everywhere. But it is equally well-known that, as soon as 
the spectrum $Sp\,( f)$ of $f$ (the set of integers $k$ at which the Fourier coefficients of $f$ do not vanish, 
{\it i.e.} $\hat f(k) \neq 0$) is sufficiently ``lacunary'', in the sense of Hadamard {\it e.g.}, then the Fourier 
series of $f$ is absolutely convergent if $f$ is continuous and almost everywhere convergent if $f$ is merely 
integrable (and in this latter case $f\in L^p$ for every $p < \infty $). Those facts have given birth to the 
theory of thin sets $\Lambda $ of integers, initiated by Rudin \cite{Ru}: those sets $\Lambda $ such that,
if $Sp\, (f ) \subseteq \Lambda $ (we shall write $f \in {\mathscr B}_\Lambda$ when $f$ is in some Banach 
function space ${\mathscr B}$ contained in $L^1 (\T)$) and $Sp\,(f) \subseteq \Lambda$), then 
$S_n(f)$, or $f$ itself, is better behaved than in the general case. 
Let us for example recall that the set $\Lambda $ is said to be: \\
- a \emph{$p$-Sidon set} $(1 \le p < 2)$ if $\hat f \in l_p$ (and not only $\hat f \in l_2$) as soon as $f$ is 
continuous and $Sp\, (f) \subseteq \Lambda $; this amounts to an ``{\it a priori} inequality''
$\Vert \hat f \Vert_p \le C\Vert f \Vert_\infty $, for each $f\in {\mathscr C}_\Lambda $; the case $p=1$ is 
the celebrated case of Sidon ($=1$-Sidon) sets;\\
- a \emph{$p$-Rider set} $(1\le p< 2)$ if we have an {\it a priori} inequality
$\Vert \hat f \Vert_p \le C\,[\![ f ]\!]$, for every trigonometric polynomial with spectrum in $\Lambda$; here  
$[\![ f ]\!]$ is the so-called Pisier norm of $f = \sum \hat f (n) e_n$, where $e_n (x) = \e^{inx}$, {\it i.e.} 
$[\![ f ]\!] = {\mathbb E}\, \Vert f_\omega \Vert_\infty $, where 
$f_\omega = \sum \varepsilon_n (\omega ) \hat f (n) e_n$, $(\varepsilon_n)$ being an {\it i.i.d.} 
sequence of centered,  $\pm 1$-valued, random variables defined on some probability space (a Rademacher 
sequence), and where ${\mathbb E}$ denotes the expectation on that space; this apparently exotic notion 
(weaker than $p$-Sidonicity) turned out to be very useful when Rider~\cite{Ri} reformulated a result of 
Drury (proved in the course of the result that the union of two Sidon set sets is a Sidon set) under the form: 
$1$-Rider sets and Sidon sets are the same (in spite of some partial results, it is not yet known whether a 
$p$-Rider set is a $p$-Sidon set: see~\cite{LR} however, for a partial result);\\
- a \emph{set of uniform convergence} (in short a $UC$-set) if the Fourier series of each 
$f\in {\mathscr C}_\Lambda $ converges uniformly, which amounts to the inequality 
$\Vert S_n (f) \Vert_\infty \le C\Vert f\Vert_\infty $, $\forall f\in {\mathscr C}_\Lambda $; 
Sidon sets are $UC$, but the converse is false;\\
- a \emph{$\Lambda (q)$-set}, $1< q < \infty $, if every $f \in L^1_\Lambda $ is in fact in $L^q$, which 
amounts to the inequality $\Vert f \Vert_q \le C_q\Vert f \Vert_1 $, $\forall f\in L^1_\Lambda $.
Sidon sets are $\Lambda (q)$ for every $q <\infty $ (and even $C_q \le C\sqrt q$); the converse is false, 
except when we require $C_q \le C\sqrt q $ (\cite{Pi});\\
- a \emph{Rosenthal set} if every $f \in L^\infty_\Lambda$ is almost everywhere equal to a continuous 
function. Sidon sets are Rosenthal, but the converse in false.
\medskip

This theory has long suffered from a severe lack of examples: those examples were always, more or less, 
sums of Hadamard sets, and in that case the banachic properties of the corresponding 
${\mathscr C}_\Lambda $-spaces were very rigid. The use of random sets (in the sense of the selectors method)
of integers has significantly changed the situation (see \cite{LQ}, and our paper \cite{LQR}). Let us recall more 
in detail the notation and setting of our previous work \cite{LQR}. The method 
of selectors consists in the following: let $(\varepsilon_k)_{k\ge 1}$ be a sequence of independent, 
$(0,1)$-valued random variables, with respective means $\delta_k$, defined on a probability space 
$\Omega $, and to which we attach the random set of integers $\Lambda = \Lambda (\omega )$,  
$\omega \in \Omega $, defined by $\Lambda (\omega ) = \{k \ge 1\,; \ \varepsilon_k (\omega ) = 1\}$.\par

The properties of $\Lambda (\omega )$ of course highly depend on the 
$\delta_k$'s, and roughly speaking the smaller the $\delta_k$'s, the better ${\mathscr C}_\Lambda $, 
$L^1_\Lambda$, \dots\,. In \cite{L}, and then, in a much deeper way, in \cite{LQR}, relying on a probabilistic 
result of J. Bourgain on ergodic means, and on a deterministic result of F.~Lust-Piquard (\cite{L-P}) on those 
ergodic means, we had randomly built new examples of sets $\Lambda $ of integers which were both: 
locally thin from the point of view of harmonic analysis (their traces on big segments $[M_n,M_{n+1}]$ of 
integers were uniformly Sidon sets); regularly distributed from the point of view of number theory, and 
therefore globally big from the point of view of Banach space theory, in that the space ${\mathscr C}_\Lambda$ 
contained an isomorphic copy of the Banach space $c_0$ of sequences vanishing at infinity. More precisely, we 
have constructed subsets $\Lambda \subseteq \N$ which are thin in the following respects: $\Lambda $ is a 
$UC$-set, a $p$-Rider set for various $p\in [1,2[$, a $\Lambda (q)$-set for every $q < \infty $, and large 
in two respects: the space ${\mathscr C}_\Lambda $  contains an isomorphic copy of $c_0$, and, most often, 
$\Lambda $ is dense in the integers equipped with the Bohr topology.\par
Now, taking $\delta_k$ bigger and bigger, we had obtained sets $\Lambda $ which were less and less thin 
($p$-Sidon for every $p>1$, $q$-Rider, but $s$-Rider for no $s <q$, $s$-Rider for every $s>q$, but not 
$q$-Rider), and, in any case $\Lambda (q)$ for every $q <\infty$, and such that ${\mathscr C}_\Lambda$ 
contains a subspace isomorphic to $c_0$. In particular, in Theorem~II.7, page 124, and 
Theorem~II.10, page 130, we take respectively $\delta_k \approx \frac{\log k}{k}$ and 
$\delta_k \approx \frac{(\log k)^\alpha}{k (\log \log k)^{\alpha +1}}$\raise 1,5pt \hbox{,} where 
$\alpha = \frac{2(p-1)}{2-p}$ is an increasing function of $p\in [1,2)$, and which becomes $\ge 1$ as $p$ 
becomes $\ge 4/3$. The case $\delta_k = \frac{1}{k}$ would correspond (randomly) to Sidon sets 
({\it i.e.} $1$-Sidon sets).\par

After the proofs of Theorem~II.7 and Theorem~II.10, we were asking two questions: \par
1) (p. 129) Our construction is very complicated and needs a second random construction of a set $E$ inside 
the random set $\Lambda $. Is it possible to give a simpler proof?\par
2) (p. 130) In Theorem~II.10, can we keep the property for the random set $\Lambda $ to be a $UC$-set, with 
high probability, when $\alpha > 1$ (equivalently when $p > \frac{4}{3}$)?
\medskip

The goal of this work is to answer affirmatively the first question (relying on a recent deviation inequality of 
Boucheron, Lugosi and Massart \cite{BLM}) and negatively the second one (relying on an older result of Kashin 
and Tzafriri \cite{KT}). This work is accordingly divided into three parts. In Section~2, we prove a (one-sided) 
concentration inequality for norms of Rademacher sums. In Section~3, we apply the concentration inequality
to get a substantially simplified proof of Theorem II.7 in \cite{LQR}. Finally, in Section~4, we give a 
(stochastically) negative answer to question~2 when $p > \frac{4}{3}$: almost surely, $\Lambda $ will not be 
a $UC$-set; here, we use the above mentionned result of Kashin and Tzafriri \cite{KT} on the non-$UC$ 
character of big random subsets of integers.

\section{\hskip -5,2 pt A one-sided inequality for norms of Rademacher sums}

Let $E$ be a (real or complex) Banach space, $v_1, \ldots, v_n$ be vectors of $E$, $X_1,\ldots , X_n$ be 
independent, real-valued, centered, random variables, and let $Z = \big\| \sum_1^n X_j v_j \|$.\par

If $\vert X_j \vert \le 1$ {\it a.s.}, it is well-known (see \cite{LT}) that:
\begin{equation}\label{eq2.1}
\P (\vert Z - \E (Z)\vert > t) \le 2 \exp 
\bigg(-\frac{t^2}{8 \sum_1^{n \phantom{\tilde l}}\Vert v_j \Vert ^2}\bigg)\,,\quad \forall t > 0.
\end{equation}

But often, the ``strong'' $l_2$-norm of the $n$-tuple $v = (v_1,\ldots, v_n)$, namely 
$\Vert v \Vert_{strong} = (\sum_{j=1}^n \Vert v_j\Vert^2)^{1/2}$, is too large for (\ref{eq2.1}) 
to be interesting, and it is advisable to work with the ``weak'' $l_2$-norm of $v$, defined by:
\begin{equation}\label{eq2.2}
\sigma = \Vert v \Vert_{weak} 
= \sup_{\varphi \in B_{E^\ast}}
\Big(\sum_1^n \vert \varphi (v_j)\vert^2\Big)^{1/2} 
= \sup_{\sum \vert a_j \vert^2 \le 1} \Big\Vert \sum^n_1 a_j v_j \Big\Vert ,
\end{equation}
where $B_{E^\ast}$ denotes the closed unit ball of the dual space $E^\ast $.\par

If $(X_j)_j$ is a standard gaussian sequence ($\E X_j = 0, \E X_j^2 = 1$), this is what Maurey and Pisier 
suceeded in doing, using either the It\^o formula or the rotational invariance of the $X_j$'s; they proved the 
following (see \cite{LQ}, Chapitre~8, Th\'eor\`eme~I.4):
\begin{equation}\label{eq2.3}
\P (\vert Z - \E Z\vert > t) \le 2\exp \Big( - \frac{t^2}{C\sigma^2}\Big)\,, \quad \forall t > 0,
\end{equation}
where $\sigma $ is as in (\ref{eq2.2}), and $C$ is a numerical constant, \textit{e.g.} $C =\pi^2/2$.\par
\smallskip

To the best of our knowledge, no inequality as simple and direct as (\ref{eq2.3}) is available for non-gaussian 
(\textit{e.g.} for Rademacher variables) variables, although several more complicated deviation 
inequalities are known: see \textit{e.g.} \cite{JS}, \cite{LT}.\par

For the applications to Harmonic analysis which we have in view, where we use the so-called 
``selectors method'', we precisely need an analogue of  (\ref{eq2.3}), in the non-gaussian, uniformly bounded 
(and centered) case; we shall prove that at least a one-sided version of (\ref{eq2.3}) holds in this case, by 
showing the following result, which is interesting for itself.

\begin{Theo}\label{theorem 2.1}
With the previous notations, assume that $\vert X_j \vert \le 1$ a.s.\,. Then, we have the one-sided estimate: 
\begin{equation}\label{eq2.4}
\P (Z - \E Z > t) \le \exp \Big( - \frac{t^2}{C\sigma^2}\Big)\,, \quad 
\forall t > 0,
\end{equation}
where $C > 0$ is a numerical constant ($C = 32$, for example).
\end{Theo}

The proof of (\ref{eq2.4}) will make use of a recent deviation inequality due to Boucheron, Lugosi and Massart 
\cite{BLM}. Before stating this inequality, we need some notation.\par 
Let $X_1,\ldots, X_n$ be independent, real-valued random variables (here, we temporarily forget the 
assumptions of the previous Theorem), and let $(X'_1,\ldots , X'_n)$ be an independent copy of 
$(X_1,\ldots , X_n)$.\par 
If $f \colon \R^n  \to \R$ is a given measurable function, we set $Z = f(X_1,\ldots , X_n)$ and 
$Z'_i = f(X_1, \ldots , X_{i-1}, X'_i, X_{i+1},\ldots , X_n)$, $1 \leq i \leq n$. 
With those notations, the Boucheron-Lugosi-Massart Theorem goes as follows:

\begin{Theo}\label{theorem 2.2}
Assume that there is some constant $a, b \geq 0$, not both zero, such that: 
\begin{equation}\label{eq2.5}
\sum_{i=1}^n (Z - Z'_i)^2 \ind_{(Z > Z'_i)} \le aZ + b  \quad \textit{a.s.}
\end{equation}
Then, we have the following one-sided deviation inequality:
\begin{equation}\label{eq2.6}
\P (Z > \E Z + t) \le \exp 
\Big(- \frac{t^2}{4a\,\E Z + 4b+2at)}\Big)\,, \quad \forall t > 0.
\end{equation}
\end{Theo}

\noindent\textbf{Proof of Theorem \ref{theorem 2.1}.}
We shall in fact use a very special case of Theorem~\ref{theorem 2.2}, the case when $a=0$; but, as the 
three fore-named authors remark, this special case is already very useful, and far from trivial to prove! 
To prove (\ref{eq2.4}), we are going to check that, for 
$f(X_1,\ldots ,X_n) = \Vert \sum_1^n X_j v_j\Vert = Z$,
the assumption (\ref{eq2.5}) holds for $a = 0$ and $b = 4\sigma^2$. 
In fact, fix $\omega \in \Omega $ and denote by $ I = I_\omega $ the set of indices $i$ such that 
$Z (\omega ) > Z'_i (\omega )$. For simplicity of notation, we assume that the Banach space $E$ is real. Let
$\varphi = \varphi_\omega \in E^\ast $ such that $\Vert \varphi \Vert = 1$ and
$Z = \varphi \big(\sum_{j=1}^n X_j v_j) = \sum_ {j=1}^n X_j \varphi (v_j)$.\par

For $i \in I$, we have 
$Z'_i(\omega ) = Z'_i \ge \varphi \big(\sum_{j\neq i} X_j v_j + X'_i v_i\big)$, so that 
$0 \le Z - Z'_i \le \sum_{j=1}^n X_j \varphi (v_j) - \sum_{j\neq i} X_j \varphi (v_j) - X'_i \varphi (v_i)
= (X_i - X'_i) \varphi (v_i)$, implying $(Z - Z'_i)^2 \le 4 \vert \varphi (v_i)\vert^2$. By summing those 
inequalities, we get: 
\begin{align*}
\sum_{i=1}^n (Z - Z'_i)^2 \ind_{(Z > Z'_i)} 
& = \sum_{i\in I} (Z - Z'_i)^2  \le 4 \sum_ {i\in I} \vert \varphi (v_i)\vert^2 
\le 4 \sum_{i=1}^n \vert \varphi (v_i)\vert^2\le 4\sigma^2 \\
& = 0.Z + 4 \sigma^2.
\end{align*}
Let us observe the crucial role of the ``conditioning'' $Z > Z'_i$ when we want to check that (\ref{eq2.5}) holds. 
Now, (\ref{eq2.4}) is an immediate consequence of (\ref{eq2.6}). \hfill $\square $
\vspace*{2mm}

\section{Construction of 4/3-Rider sets}

We first recall some notations of \cite{LQR}. $\Psi_2$ denotes the Orlicz function $\Psi_2 (x) = \e^{x^2}-1$, 
and $\Vert \; \Vert_{\Psi_2}$ is the corresponding Luxemburg norm. If $A$ is a finite subset of the 
integers, $\Psi_A$ denotes the quantity $\Vert \sum_{n\in A} e_n \Vert_{\Psi_2}$, 
where $e_n(t) = \e^{int}$, $t\in \R/2\pi  \Z = \T $, and $\T$ is equipped with its Haar measure $m$. 
$\Lambda $ will always be a subset of the positive integers $\N$. Recall that $\Lambda $ is 
\emph{uniformly distributed} if the ergodic means 
$A_N(t) = \frac{1}{\vert\Lambda_N\vert}\sum_{n\in \Lambda _N} e_n(t)$ tend to zero as $N \to \infty $,
for each $t\in \T$, $t\neq 0$. Here, $\Lambda_N = \Lambda \cap [1,N]$.
If $\Lambda $ is uniformly distributed, ${\mathscr C}_\Lambda $ contains $c_0$, and if 
${\mathscr C}_\Lambda $ contains $c_0$, $\Lambda $ cannot be a Rosenthal set (see \cite{LQR}).
According to results of J. Bourgain (see \cite{LQR}) and F. Lust-Piquard (\cite{L-P}), respectively, a random set 
$\Lambda $ corresponding to selectors of mean $\delta_k$ with $k\delta_k \to  \infty $ is almost surely 
uniformly distributed and if a subset $E$ of a uniformly distributed set $\Lambda $ has positive upper density 
in $\Lambda $, {\it i.e.} if 
$\limsup_N \frac{\vert E \cap [1,N]\vert}{\vert \Lambda \cap [1,N]}  > 0$,
then ${\mathscr C}_E$ contains $c_0$, and $E$ is non-Rosenthal.
\smallskip

In  \cite{LQR}, we had given a fairly complicated proof of the following theorem (labelled as Theorem~II.7):

\begin{Theo}\label{theorem 3.1}
There exists a subset $\Lambda $ of the integers, which is uniformly distributed, and contains a subset 
$E$ of positive integers with the following properties:\par
1) $E$ is a $\frac{4}{3}$-Rider set, but is not $q$-Rider for $q <4/3$, a $UC$-set, and a $\Lambda (q)$-set 
for all $q < \infty $; \par
2) $E$ is of positive upper density inside $\Lambda $; in particular, ${\mathscr C}_E$ contains $c_0$ and 
$E$ is not a Rosenthal set.
\end{Theo}

We shall show here that the use of Theorem \ref{theorem 2.1} allows a substantially simplified proof, which 
avoids a double random selection. We first need the following simple lemma.

\begin{Lem}\label{lemma 3.2}
Let $A$ be a finite subset of the integers, of cardinality $n \ge 2$; let $v = (e_j)_{j\in A}$, considered as an 
$n$-tuple of elements of the Banach space $E = L^{\Psi_2} = L^{\Psi_2}(\T , m)$, and let $\sigma$
be its weak $l_2$-norm. Then:
\begin{equation}\label{eq3.1}
\smash{\sigma \le C_0 \sqrt{ \frac{n}{\log n}  }}\,,
\end{equation}
where $C_0$ is a numerical constant.
\end{Lem}

\noindent\textbf{Proof.} 
Let $a = (a_j)_{j\in A}$ be such that $\sum_{j\in A} \vert a_j\vert^2 = 1$. Let 
$f = f_a = \sum_{j\in A} a_j e_j$, and $M = \Vert f\Vert_\infty $. By H\"older's inequality, we have 
$\frac{ \Vert f \Vert_p}{\sqrt p}  \le \frac{M}{\sqrt{p}\,M^{2/p}}$ for $2 < p < \infty $. Since 
$M \leq \sqrt n$, we get
$\frac{\Vert f\Vert_p}{\sqrt p} \le  \frac{\sqrt n}{\sqrt{p}\,n^{1/p}} \le C \sqrt{\frac{n}{\log n}}$. 
By  Stirling's formula, $\Vert f \Vert_{\Psi_2} \approx \sup_{p>2} \frac{ \Vert f \Vert_p}{\sqrt p}$, so 
the lemma is proved, since $\sigma= \sup_a \| f_a\|_{\Psi_2}$\hfill $\square $
\medskip

We now turn to the shortened proof of Theorem~\ref{theorem 3.1}.
\medskip

Let $I_n = [2^n,2^{n+1}[$, $n \ge 2$ ; $\delta_k = c \frac{n}{2^n}$ if $k\in I_n$ $(c > 0)$.\par\smallskip
Let $(\varepsilon _k)_k$ be a sequence of ``selectors'', \textit{i.e.} independent, $(0,1)$-valued, random 
variables of expectation $\E\varepsilon_k = \delta_k$, and let $\Lambda = \Lambda (\omega )$ be the 
random set of positive integers defined by $\Lambda = \{ k \ge 1\,;\  \varepsilon_k = 1 \}$.
We set also $\Lambda_n = \Lambda \cap I_n$ and  
$\sigma_n = \E \vert \Lambda_n \vert =\sum_{k\in I_n} \delta_k = cn$.\par\smallskip

We shall now need the following lemma (the notation $\Psi_A$ is defined at the beginning of the section).

\begin{Lem}\label{lemma 3.3}
Almost surely, for $n$ large enough: 
\begin{align}
\frac{c}{2} n \le \vert \Lambda_n \vert \le 2cn \ ;\label{eq3.4}
\\
\Psi_{\Lambda_n}\le  C''\vert \Lambda_n \vert ^{1/2} \ .\label{eq3.5}
\end{align}
\end{Lem}

\noindent\textbf{Proof :}
(\ref{eq3.4}) is the easier part of Lemma~II.9 in \cite{LQR}. To prove (\ref{eq3.5}), we recall an inequality due 
to G. Pisier \cite{Pi}: if $(X_k)$ is a sequence of independent, centered and square-integrable, random variables 
of respective variances $V(X_k)$, we have:
\begin{equation}\label{eq3.6}
\E\Big\Vert \sum_k X_k e_k \Big\Vert_{\Psi_2}  \le C_1 \Big(\sum_k V(X_k)\Big)^{1/2}.
\end{equation}
Applying (\ref{eq3.6}) to the centered variables $X_k = \varepsilon_k - \delta_k$, we get, assuming $c\le 1$: 
\begin{displaymath}
\E\Big\Vert \sum_{k\in I_n}(\varepsilon_k - \delta_k) e_k \Big\Vert_{\Psi_2}
\le C_1 \Big(\sum_{k\in I_n} \delta_k(1 - \delta_k) \Big)^{1/2}
\le C_1 \Big(\sum_{k\in I_n} \delta_k \Big)^{1/2}\le C_1 \sqrt n .
\end{displaymath}
Now, set $Z_n = \big\Vert \sum_{k\in I_n} (\varepsilon_k -\delta_k) e_k \big\Vert_{\Psi_2}$. Let 
$\lambda $ be a fixed real number $> 1$, and $C_0$ be as in Lemma~\ref{lemma 3.2}. Applying 
Theorem~\ref{theorem 2.1} with $C = 32$, and $t_n = \lambda\sqrt{32 C_0^2n}$, we get, using 
Lemma~\ref{lemma 3.2}:
\begin{displaymath}
\P (Z_n - \E Z_n > t_n) \le \exp \Big(- \frac{t^2_n}{32 \sigma^2}\Big)
\le \exp \Big(- \frac{32\lambda^2 C_0^2 n \log n}{32 C_0^2 n}\Big)
= n^{-\lambda^2}.
\end{displaymath}
By the Borel-Cantelli Lemma, we have almost surely, for $n$ large enough:
\begin{displaymath}
Z_n \le \E Z_n + t_n \le (C_1+4C_0\lambda ) \sqrt n = C_2 \sqrt n.
\end{displaymath}
For such $\omega $'s and $n$'s, it follows that: 
\begin{align*}
\Psi_{\Lambda_n} = \Big\Vert \sum_{k\in I_n}\varepsilon_k e_k \Big\Vert_{\Psi_2} 
&\le Z_n + \Big\Vert \sum_{k\in I_n} \delta_k e_k \Big\Vert_{\Psi_2}
\le Z_n + \frac{n}{2^n}\Big\Vert \sum_{k\in I_n} e_k\Big\Vert_{\Psi_2}\\
&\le C_2 \sqrt n +  \frac{n}{2^n} C_0 \frac{2^n}{\sqrt {\log 2^n}}
= : C_3 \sqrt n,
\end{align*}
because, with the notations of Lemma~\ref{lemma 3.2}, we have:
\begin{displaymath}
\Big\Vert \sum_{k\in I_n} e_k \Big\Vert_{\Psi_2}
\le \sqrt{\vert I_n\vert } \sigma \le 2^{n/2} C_0 \frac{ 2^{\frac{n}{2}}}{\sqrt {\log 2^n}} \cdot 
\end{displaymath}
This ends the proof of Lemma~\ref{lemma 3.3}, because we know that
$n\le \frac{2}{c}\vert\Lambda_n\vert $ for large $n$, almost surely, and therefore
$\Psi_{\Lambda_n} \le C_3 \sqrt{\frac{2}{c}} \vert \Lambda _n\vert^{1/2}
=: c''\vert \Lambda_n\vert^{1/2}$, {\it a.s.}\,. \hfill $\square $
\medskip

We now prove Theorem~\ref{theorem 3.1} as follows: let us fix a point $\omega\in \Omega $ in such a way 
that $\Lambda = \Lambda (\omega )$ is uniformly distributed and that $\Lambda_n $ verifies (\ref{eq3.4}) 
and  (\ref{eq3.5}) for $n \ge n_0$; this is possible from \cite{LQR} and from Lemma \ref{lemma 3.3}. 
We then use a result of the third-named author (\cite{Ro}), asserting that there is a numerical constant 
$\delta >0$ such that each finite subset $A$ of $\Z^\ast $ contains a quasi-independent subset $B$ such that 
$\vert B\vert \ge \delta \big( \frac{\vert A\vert}{\,\Psi_A}\big)^2$ (recall that a subset $Q$ of $\Z$ is said to 
be quasi-independent if, whenever $n_1, \ldots, n_k \in Q$, the equality $\sum_{j=1}^k \theta_j n_j =0$ with 
$\theta_j =0, -1, +1$ holds only when $\theta_j=0$ for all $j$). This allows us to select inside each 
$\Lambda_n $ a quasi-independent subset $E_n$ such that:
\begin{equation}\label{eq3.7}
\vert E_n\vert \ge \delta 
\Big(\frac{\vert \Lambda_n\vert}{\Psi _{\Lambda_n}}\Big)^2 \ge
 \frac{\delta}{{c''}^2} \vert \Lambda_n \vert = : \delta '\vert\Lambda_n\vert\,.
\end{equation}
A combinatorial argument (see \cite{LQR}, p. 128--129) shows that, if $E = \cup_{n > n_0} E_n$, 
then each finite $ A\subset E$ contains  a quasi-independent subset $B\subseteq A$ such that  
$\vert B\vert \ge \delta \vert A\vert^{1/2}$. By \cite{Ro}, $E$ is a $\frac{4}{3}$-Rider set. The set $E$ has 
all the required properties. Indeed, it follows from Lemma~\ref{lemma 3.2},~a) that 
$\vert E \cap [1,N] \vert \ge \delta (\log N)^2$. If now $E$ is $p$-Rider, we must have 
$\vert E \cap [1,N]\vert \le C (\log N)^{\frac{p}{2-p}}$; therefore $2 \leq \frac{p}{2-p}$\raise 1,5pt \hbox{,} 
so $p \ge 4/3$. The fact that $E$ is both $UC$ and $\Lambda (q)$ is due to the local character of these 
notions, and to the fact that the sets $E \cap [2^n,2^{n+1}[ = E_n$ are by construction quasi-independent 
(as detailed in \cite{LQR}). On the other hand, since each $E_n$ is approximately proportional to 
$\Lambda_n$, $E$ is of positive upper density in $\Lambda $. Now $\Lambda $ is uniformly distributed 
(by Bourgain's criterion: see \cite{LQR}, p.~115). Therefore, by the result of F. Lust-Piquard (\cite{L-P}, and 
see Theorem~I.9, p.~114 in \cite{LQR}), ${\mathscr C}_E$ contains $c_0$, which prevents $E$ from being 
a Rosenthal set. \hfill $\square $

\section {$p$-Rider sets, with $p > 4/3$, which are not  $UC$-sets}

Let $p \in ]\frac{4}{3}, 2[$, so that $\alpha = \frac{2(p-1)}{2-p} > 1$. As we mentioned in the Introduction, 
the random set $\Lambda = \Lambda (\omega )$ of integers in Theorem~II.10 of \cite{LQR} corresponds to 
selectors $\varepsilon_k$ with mean
$\delta_k = c \frac{(\log k)^\alpha}{k (\log \log k)^{\alpha +1}}\,\cdot$ We shall prove the following:

\begin{Theo}\label{theorem 4.1}
The random set $\Lambda $ corresponding to selectors of mean 
$\delta_k = c \frac{(\log k)^\alpha}{k (\log \log k)^{\alpha +1}}$ has almost surely the following properties:
\par
a) $\Lambda $ is $p$-Rider, but $q$-Rider for no $q < p$;\par
b) $\Lambda $ is $\Lambda (q)$ for all $q < \infty $;\par
c) $\Lambda $ is uniformly distributed; in particular, it is dense in the Bohr group and 
${\mathscr C}_\Lambda $ contains $c_0$;\par
d) $\Lambda $ is {\bf not} a $UC$-set.
\end{Theo}

\noindent\textbf{Remark.} This supports the conjecture that $p$-Rider sets with $p > 4/3$ are not
of the same nature as $p$-Rider sets for $p < 4/3$ (see also \cite{LLQR}, Theorem~3.1. and \cite{LR}). \par
\smallskip

The novelty here is \textit{d)}, which answers in the negative a question of \cite{LQR} and we shall mainly 
concentrate on it, although we shall add some details for \textit{a)},\textit{b)}, \textit{c)}, since the proof
 of Theorem~II.10 in \cite{LQR} is too sketchy and contains two small misprints (namely ($\ast $) and 
($\ast \ast $), p.~130).
\smallskip

Recall that the $UC$-constant $U(E)$ of a set $E$ of positive integers is the smallest constant $M$ such that 
$\Vert S_N f\Vert_\infty \le M\Vert f\Vert_\infty $ for every $f\in {\mathscr C}_E$ and every non-negative 
integer $N$, where $S_N f = \sum_{-N}^N\hat f(k) e_k$. We shall use the following result of 
Kashin and Tzafriri \cite{KT}:
\begin{Theo}\label{theorem 4.2}
Let $N \ge 1$ be an integer and $\varepsilon '_1 ,\ldots, \varepsilon '_N$ be selectors of equal mean $\delta $. 
Set $\sigma (\omega ) = \{k\le N \,;\  \varepsilon'_k (\omega ) = 1\}$. Then:
\begin{equation}\label{eq4.2}
\smash{\P \Big(U\big(\sigma (\omega )\big) 
\le \gamma \log \Big(2 + \frac{\delta N}{\log N} \Big) \Big)
\le \frac{5}{N^3} \raise 1pt \hbox{,}}
\end{equation}
where $\gamma $ is a positive numerical constant.
\end{Theo}

We now turn to the proof of Theorem~\ref{theorem 4.1}. As in \cite{LQR}, we set, for a fixed $\beta > \alpha$:
\begin{equation}\label{eq4.3}
\smash{M_n = n^{\beta n}\;; \quad
\Lambda_n = \Lambda \cap [1,n] \; ; \quad  
\Lambda^\ast_n = \Lambda \cap [M_n,M_{n+1}[. }
\end{equation}
We need the following technical lemma, whose proof is postponed
(and is needed only for \textit{a)},  \textit{b)}, \textit{c)}).
\begin{Lem}\label{lemma 4.3}
We have almost surely for large $n$
\begin{equation}\label{eq4.4}
\vert \Lambda_{M_n}\vert \approx n^{\alpha +1}\; ; \qquad
\vert \Lambda^\ast_n \vert \approx n^\alpha .
\end{equation}
\end{Lem}

Observe that, for $k \in \Lambda^\ast_n$, one has:
\begin{displaymath}
\delta_k = 
c\,\frac{(\log k)^\alpha}{k(\log\log k)^{\alpha+1}} \gg
\frac{ (n \log n)^\alpha}{M_{n+1}(\log n)^{\alpha +1}}
= \frac{n^\alpha}{M_{n+1}\log n} =: \frac{q_n}{N_n}\,\raise 1pt \hbox{,}
\end{displaymath}
where $N_n = M_{n+1}- M_n$ is the number of elements of the support of $\Lambda^\ast_n$ (note that 
$N_n \sim M_{n+1}$), and where $q_n$ is such that
\begin{equation}\label{eq4.6}
\smash{q_n \approx \frac{\,n^\alpha}{\log n}\,\cdot }
\end{equation}
We can adjust the constants so as to have $\delta _k \ge q_n/N_n$ for $k \in \Lambda^\ast_n$. Now, we 
introduce selectors $(\varepsilon ''_k)$ independent of the $\varepsilon_j$'s, of respective means 
$\delta''_k = q_n/(N_n \delta_k)$. Then the selectors $\varepsilon'_k = \varepsilon_k \varepsilon''_k$ 
have means $\delta '_k = q_n/N_n$ for $k\in \Lambda^\ast_n$, and we have 
$\delta_k \ge \delta '_k$ for each $k \ge 1$.\par

Let $\Lambda ' = \{k\,;\  \varepsilon '_k = 1\}$ and $ {\Lambda '}_n^\ast  = \Lambda '\cap [M_n,M_{n+1}[$. 
It follows from (\ref{eq4.2}) and the fact that $U(E+a) = U(E)$ for any set $E$ of positive integers and any 
non-negative integer $a$ that:
\begin{displaymath}
\smash{\P \Big(U({\Lambda '}^\ast_n) \le \gamma \log 
\Big(2 + \frac{q_n}{\log N_n}\Big)\Big) \le 5N_n^{-3}.}
\end{displaymath}
By the Borel-Cantelli Lemma, we have almost surely 
$U({\Lambda '}^\ast _n) > \gamma \log \big(2 + \frac{q_n}{\log N_n}\big)$ for $n$ large enough. 
But we see from (\ref{eq4.4}) and (\ref{eq4.3}) that:
\begin{displaymath}
\smash{\frac{q_n}{\log N_n}\approx
\frac{n^\alpha}{(\log n) (n\log n)} = \frac{n^{\alpha -1}}{(\log n)^2}\, \raise 1,5pt \hbox{,}}
\end{displaymath}
and this tends to infinity since $\alpha > 1$. This shows that $\Lambda '$ is almost surely non-$UC$. And due 
to the construction of the $\varepsilon'_k$'s, we have: $\Lambda \supseteq \Lambda '$ almost surely. 
This of course implies that $\Lambda $ is not a $UC$-set either (almost surely), ending the proof 
of \textit{d)} in Theorem~\ref{theorem 4.1}.\hfill $\square$
\medskip

We now indicate a proof of the lemma. Almost surely, $\vert \Lambda_{M_n}\vert $ behaves for large $n$ as: 
\begin{align*}
\E (\vert \Lambda_{M_n}\vert )
&= \sum^{M_n}_1 
\frac{(\log k)^\alpha}{k(\log\log k)^{\alpha +1}} 
\approx \int^{M_n}_{e^2} \frac{(\log t)^\alpha}{t(\log\log t)^{\alpha +1}}dt\\
& = \int^{\log M_n}_2 \hskip -8pt  \frac{x^\alpha dx}{(\log x)^{\alpha +1}}
\approx \frac{1}{(\log n)^{\alpha +1}} \int^{\log M_n}_2 \hskip -15pt x^\alpha dx 
\approx \frac{(\log M_n)^{\alpha +1}}{(\log n)^{\alpha +1}} \approx n^{\alpha +1}.
\end{align*}
Similarly, $\vert \Lambda^\ast_n \vert $ behaves almost surely as:
\begin{align}
\int^{M_{n+1}}_{M_n} \hskip - 8pt 
\frac{(\log t)^\alpha}{t(\log\log t)^{\alpha+1}} dt 
& = \int_{\log M_n}^{\log M_{n+1}} \frac{x^\alpha}{(\log x)^{\alpha +1}} dx 
\approx \frac{1}{(\log n)^{\alpha +1}} x^\alpha dx \notag \\ 
&\approx  \frac{1}{(\log n)^{\alpha +1}}(\log M_{n+1} - \log M_n) (\log
M_n)^\alpha  \notag \\ 
& \approx \frac{1}{(\log n)^{\alpha +1}}\log n (n \log n)^\alpha
 \approx n^\alpha . \tag*{$\square$}
\end{align} 
To finish the proof, we shall use a lemma of \cite{LQR} (recall that a \emph{relation of length $n$} in 
$A \subseteq \Z^\ast$ is a $(-1, 0, +1)$-valued sequence $(\theta_k)_{k\in A}$ such that 
$\sum_{k\in A} \theta_k\,k =0$ and  $\sum_{k\in A} |\theta_k| =n$):
\begin{Lem}\label{lemma 4.2}
Let $n\ge 2$ and $M$ be integers. Set 
\begin{displaymath}
\smash{ \Omega_n (M) = \{\omega \mid \Lambda (\omega) \cap [M,\infty [\,
\textrm{contains at least a relation of length }  n\}.}
\end{displaymath}
Then:
\begin{displaymath}
\smash{\P [\Omega_n(M)] \le \frac{C^n}{n^n}\sum_{j> M}\delta^2_j \sigma_j^{n-2},}
\end{displaymath}
\vskip 5pt 
\noindent where $\sigma_j = \delta_1 + \ldots + \delta_j$, and $C$ is a numerical constant.
\end{Lem}

In our case, with $M = M_n$, this lemma gives :
\begin{align*}
\P [\Omega_n(M)] 
&\ll \frac{C^n}{n^n} \sum_ {j > M}
\frac{(\log j)^{2\alpha}}{j^2 (\log\log j)^{2\alpha + 2}}
\left[\frac{(\log j)^{\alpha +1}}{(\log \log j)^{\alpha + 1}}\right]^{n-2}\\
&\ll  \frac{C^n}{n^n} \int^\infty_M
\frac{(\log t)^{(\alpha +1)n+2\alpha}}
{(\log \log t)^{(\alpha+1)n+2\alpha +2}}
\frac{dt}{t^2}
\end{align*}
and an integration by parts (see \cite{LQR}, p.~117--118) now gives: 
\begin{align*}
\P [\Omega_n(M)] 
&\ll \frac{C^n}{n^n} \frac{1}{M}
\frac{(\log M)^{(\alpha +1)n+2\alpha}} {(\log\log M)^{(\alpha+1)n+2\alpha+2}}\\
&\ll \frac{C^n}{n^n} \frac{1}{n^{\beta n}} 
\frac{(n\log n)^{(\alpha+1)n +2\alpha}} {(\log n)^{(\alpha +1)n+2\alpha +2}}
\ll \frac{n^{2\alpha} C^n}{n^{(\beta -\alpha )n }(\log n)^2} ;
\end{align*}
then the assumption $\beta > \alpha $ (which reveals its importance here!) shows that 
$\sum_n  \P [\Omega_n(M_n)] < \infty $, so that, almost surely
$\Lambda (\omega) \cap [M_n,\infty [$ contains no relation of length $n$, for $n \ge n_0$. 
Having this property at our disposal, we prove (exactly as in \cite{LQR}, p.~119--120) that $\Lambda $ is 
$p$-Rider. It is not $q$-Rider for $q < p$, because then 
$\vert \Lambda_{M_n}\vert \ll (\log M_n)^{\frac{q}{2-q}} \ll (n\log n)^{\frac{q}{2-q}},$
whereas  (\ref{eq4.4}) of Lemma \ref{lemma 4.3} shows that
$\vert \Lambda_{M_n}\vert \gg n^{\alpha +1}$, with $\alpha +1 = \frac{p}{2-p} > \frac{q}{2-q}\cdot$
This proves \textit{a}). Conditions \textit{b)},\textit{c)} are clearly explained in \cite{LQR}. \hfill $\square $


\noindent {\scriptsize  Daniel  Li:}{\textit{\scriptsize 
Université d'Artois, 
Laboratoire de Math\'ematiques de Lens EA 2462--FR 2956, Facult\'e des Sciences Jean Perrin, 
23, rue J. Souvraz SP 18, F-62307 Lens Cedex, France --\\ 
daniel.li@euler.univ-artois.fr}}\\
{\scriptsize  Herv\'e Queff\'elec: } {\textit{\scriptsize 
Laboratoire Paul Painlev\'e UMR CNRS 8524, U.F.R. de Math\'ematiques Pures et Appliqu\'ees, B\^at. M2, 
Universit\'e des Sciences et Technologies de Lille 1, F-59665 Villeneuve d'Ascq Cedex, France -- 
Herve.Queffelec@math.univ-lille1.fr}}\\
{\scriptsize  Luis Rodr{\'\i}guez-Piazza:} {\textit{\scriptsize
Universidad de Sevilla, Facultad de Matematicas, Departamento de An\'alisis Matem\'atico,
Apartado de Correos 1160, 41080 Sevilla, Spain --  piazza@us.es}}

\end{document}